\documentclass[10pt]{article}

\usepackage{graphics}
\usepackage{epsfig}
\usepackage{amssymb}
\usepackage{amsthm}
\usepackage{mathrsfs}
\usepackage{hhline}
\usepackage{amsmath}
\usepackage{array}
\usepackage{float}
\usepackage{picinpar}
\usepackage{multirow}
\usepackage{cite}
\usepackage{enumerate}
\usepackage{xcolor}

\newcommand{\mb}[1]{\text{\boldmath ${#1}$}}
\def\blue{\color{blue}}

\newtheorem{theorem}{Theorem}[section]
\newtheorem{lemma}[theorem]{Lemma}
\newtheorem{example}{Example}[section]
\newtheorem{corollary}[theorem]{Corollary}
\newtheorem{problem}{Problem}[section]

\newtheorem{conjecture}{Conjecture}[section]

\newtheorem{defi}{Definition}[section]

\newenvironment{definition}{\begin{defi}\rm }{\end{defi}}

\usepackage{graphicx}
\graphicspath{%
    {converted_graphics/}
    {/}
}
\begin{document}
\baselineskip18truept
\normalsize
\begin{center}
{\mathversion{bold}\Large \bf Further Results On $k$-Super Graceful Graphs}

\bigskip
{\large Gee-Choon Lau\footnote{Corresponding author. E-mail: geeclau@yahoo.com}}\\
\emph{Faculty of Computer \& Mathematical Sciences,}\\
\emph{Universiti Teknologi MARA Johor,}\\
\emph{85009, Segamat, Malaysia.}
\medskip

{\large Wai-Chee Shiu\footnote{E-mail: wcshiu@associate.hkbu.edu.hk}}\\
\emph{Department of Mathematics, The Chinese University of Hong Kong, Shatin, Hong Kong,}\\
\emph{College of Global Talents, Beijing Institute of Technology, Zhuhai, China.}

\medskip

{\large Ho-Kuen Ng\footnote{E-mail: ho-kuen.ng@sjsu.edu}}\\
\emph{Department of Mathematics, San Jos\'{e} State University,}\\
\emph{San Jose CA 95192 USA.}

\medskip

{\large Zhen-Bin Gao\footnote{E-mail: gaozhenbin@aliyun.com}}\\
\emph{College of Mathematical Sciences, Harbin Engineering University,}\\
\emph{Harbin, 150001, P.R. China.}

\medskip

{\large Karl Schaffer\footnote{E-mail: karl\_schaffer@fhda.edu}}\\
\emph{Department of Mathematics, De Anza College,}\\
\emph{Cupertino, CA95014, USA.}

\end{center}


\medskip
\begin{abstract}
Let $G=(V(G),E(G))$ be a simple, finite and undirected graph of order $p$ and size $q$. For $k\ge 1$, a bijection $f: V(G)\cup E(G) \to  \{k, k+1, k+2,  \ldots, k+p+q-1\}$ such that $f(uv)= |f(u) - f(v)|$ for every edge $uv\in E(G)$ is said to be a $k$-super graceful labeling of $G$. We say $G$ is $k$-super graceful if it admits a $k$-super graceful labeling. In this paper, we study the $k$-super gracefulness of some graphs in which each component is either regular or bi-regular. 

\medskip

\end{abstract}

\tolerance=10000
\baselineskip12truept
\def\qed{\hspace*{\fill}$\Box$\medskip}

\def\s{\,\,\,}
\def\ss{\smallskip}
\def\ms{\medskip}
\def\bs{\bigskip}
\def\c{\centerline}
\def\nt{\noindent}
\def\ul{\underline}
\def\lc{\lceil}
\def\rc{\rceil}
\def\lf{\lfloor}
\def\rf{\rfloor}
\def\a{\alpha}
\def\b{\beta}
\def\n{\nu}
\def\o{\omega}
\def\ov{\over}
\def\m{\mu}
\def\t{\tau}
\def\th{\theta}
\def\k{\kappa}
\def\l{\lambda}
\def\L{\Lambda}
\def\g{\gamma}
\def\d{\delta}
\def\D{\Delta}
\def\e{\epsilon}
\def\lg{\langle}
\def\rg{\tongle}
\def\p{\prime}
\def\sg{\sigma}
\def\to{\rightarrow}

\newcommand{\K}{K\lower0.2cm\hbox{4}\ }
\newcommand{\cl}{\centerline}
\newcommand{\om}{\omega}
\newcommand{\ben}{\begin{enumerate}}

\newcommand{\een}{\end{enumerate}}
\newcommand{\bit}{\begin{itemize}}
\newcommand{\eit}{\end{itemize}}
\newcommand{\bea}{\begin{eqnarray*}}
\newcommand{\eea}{\end{eqnarray*}}
\newcommand{\bear}{\begin{eqnarray}}
\newcommand{\eear}{\end{eqnarray}}







\section{Introduction}

\nt Let $G=(V(G),E(G))$ (or $G=(V,E)$ for short) be a simple, finite and undirected graph without isolated vertex of order $|V|=p$ and size $|E| = q$ $(q\ge 1)$. $G$ is also called a $(p,q)$-graph. For integers $a$ and $b$ with $a\le b$, let $\mb{[a,b]}$ be the set of integers between $a$ and $b$ inclusively. All notation not defined in this paper can be found in~\cite{Bondy}. Rosa~\cite{Rosa} defined a graceful labeling of $G$ as an injective vertex labeling function $f: V\to \mb{[0,q]}$ such that the induced edge-labeling function $f^*(uv)=|f(u)-f(v)|$ for every $uv\in E$ is also injective. The following conjecture has since then become one of the most famous unsolved graph labeling problems.

\ms\nt {\bf Conjecture.} {\it All trees are graceful.}

\ms\nt Since then, there have been more than 1500 research papers on graph labelings (see the dynamic survey by Gallian~\cite{Gallian}).

\ms\nt In~\cite{Bange}, Bange et al.  defined a $k$-sequentially additive labeling $f$ of a graph $G$ as a bijection from $V\cup E$ to $\mb{[k,k+p+q-1]}$ such that for each edge $uv\in E$, $f(uv)=f(u)+f(v)$. A graph $G$ admitting a $k$-sequentially additive labeling is called a {\it $k$-sequentially additive graph.} If $k=1$, then $G$ is called a {\it simply sequentially additive graph} or an {\it SSA-graph}. They conjectured that all trees are SSA-graphs. More results on $k$-sequentially additive labeling can be found in~\cite{Hajnal+Nagy, Hegde+Miller, Manimekalai}.

\ms\nt In~\cite{Lau+Shiu+Ng-kSG}, the concept of $k$-super graceful labeling was introduced. This was referred to as a $k$-sequential labeling in~\cite{Slater} that we are only aware of after the completion of this paper.

\begin{definition} Given $k\ge 1$, a bijection $f: V\cup E \to  \mb{[k, k+p+q-1]}$ is called a {\it $k$-super graceful labeling} if $f(uv) = |f(u)-f(v)|$ for every edge $uv$ in $G$. We say $G$ is {\it $k$-super graceful} if it admits a $k$-super graceful labeling.
\end{definition}


\ms\nt This is a generalization of a super graceful labeling defined in~\cite{Perumal1, Perumal2}. For simplicity, 1-super graceful is also known as super graceful. Properties of $k$-super gracefulness and sufficient conditions on the existence of $k$-super graceful labeling of various bipartite and tripartite graphs have been investigated~\cite{Lau+Shiu+Ng-kSG}.  In this paper, we study the $k$-super gracefulness of some graphs in which each component is either regular or bi-regular.

\section{Basic Properties}

\nt Let $G+H$ be the disjoint union of graphs $G$ and $H$. Let $nG$ be the disjoint union of $n\ge2$ copies of $G$. Suppose a graph $G$ admits a $k$-super graceful labeling with the edge-label set $A$. We will say that $G$ is $k$-super graceful with edge-label set $A$ or with edges in the set $A$. The following 3 results were obtained in~\cite{Lau+Shiu+Ng-kSG}.

\begin{theorem}\label{thm-NAV} For $k\ge 1$, if a $(p,q)$-graph $G$ admits a $k$-super graceful labeling, then the $k$ largest integers in $\mb{[k,k+p+q-1]}$ must be vertex labels of $k$ mutually non-adjacent vertices. Moreover, no two of the $k+1$ smallest integers are vertex labels of adjacent vertices.
\end{theorem}

\begin{corollary}\label{cor-NAV} If $G$ is $k$-super graceful, then $1\le k\le\alpha$, where $\alpha$ is the independence number of $G$. Moreover, the upper bound is sharp. \end{corollary}

\begin{lemma}~\label{lem-G+H} Suppose $H$ is a graph of size $q$.  If $H$ is super-graceful with edge-label set $\mb{[1,q]}$ and $G$ is $(q+1)$-super graceful, then $G+H$ is super graceful. \end{lemma}

\nt By an argument similar to the proof of Lemma~\ref{lem-G+H}, we have

\begin{lemma}\label{lem-G+H3}  If there is a sequence of graphs $G_i, i=1,2,\dots,t, t\ge 2$, such that $G_i$ is $k_i$-super graceful with largest label $r_i$ and $k_{i+1} = r_i + 1$ for $1\le i\le t-1$, then $\sum^t_{i=j} G_i$ is $k_j$-super graceful for $j\ge 1$.  \end{lemma}

\begin{lemma}\label{lem-G+H4} If there is a sequence of graphs $G_i$, $1\le i\le t$, $t\ge 2$, of sizes $q_i$ such that $G_i$ is $k_i$-super graceful with edge labels being the smallest $q_i$ labels and that $k_{i+1}=k_i+q_i$, $1\le i\le t-1$, then for $s\le t$, $\sum^s_{i=j} G_i$ is $k_j$-super graceful with edge-label set $\mb{[k_j,k_j+\sum^s_{i=j}q_i-1]}$.  Conversely, for $t\ge 2$, let $G_i$ be  a graph of size $e_i$, $1\le i\le t$. For each $k\ge 1$, if there exists a $k$-super graceful $G=\sum^t_{i=1} G_i$ of size $e$ with edge-label set $\mb{[k,k+e-1]}$ and vertex labels of each $G_i$ form a sequence of consecutive integers, then each $G_i$ is $k_i$-super graceful for some $k_i$ with edge-label set $\mb{[k_i,k_i+e_i-1]}$.  \end{lemma}

\begin{lemma}\label{lem-G+H5} Suppose $H$ is a $(p,q)$-graph that admits an injective vertex labeling $f: V(H)\to \mb{[r,r+p-1]}$, $r\ge 1$, such that $f(uv)=|f(u)-f(v)|$ for each $uv\in E(H)$ with edge-label set $\mb{[k,k+q-1]}$, then $H$ is $k$-super graceful. Moreover, if $G$ is $(k+q)$-super graceful, then $G+H$ is $k$-super graceful.    \end{lemma}

\ms\nt For connected graphs, the path $P_{n+1}$ and the star $K(1,n)$ of size $n\ge 1$ are super graceful with edge-label set $\mb{[1,n]}$. For disconnected graphs, we may apply Lemma~\ref{lem-G+H4} to construct $k$-super graceful graphs of size $q$ whose edges are labeled by the smallest $q$ labels.

\begin{example}\label{eg-H+H}{\rm  It is known that $K_2$ is super graceful with edge label 1 and that $4K_2$ is 2-super graceful with edge labels 2 to 5 and corresponding end-vertex labels $(7,9)$, $(10,13)$, $(8,12)$, $(6,11)$. By Lemma~\ref{lem-G+H4}, $5K_2$ is super graceful with edge labels 1 to 5.  Also, $P_3$ is super graceful with edge labels 1, 2 and $5K_2$ is 3-super graceful with edge labels 3 to 7 and corresponding end-vertex labels $(11, 14)$, $(9, 13)$, $(12,17)$, $(10,16)$, $(8,15)$. Hence, $P_3 + 5K_2$ is super graceful with edge labels 1 to 7. }
\end{example}

\nt Now suppose that a $(p,q)$-graph $G$ admits a $k$-super graceful labeling such that all edges of $G$ are labeled by even integers. Suppose there are $a$ odd vertex labels and $b\ge 0$ even vertex labels. Thus $p=a+b$. By the assumption mentioned at the beginning of this paper, $a\ge 2$ and $b\ne 1$. Suppose $G$ has $\omega$ components. Since there is no isolated vertex, $q\ge p-\omega \ge p-q$. Hence $2q\ge p$. Suppose $2q=p$, then $\omega =q$. This implies that $G\cong qK_2$. We will study this case in the next section.

\ms\nt Now we assume that $a+b=p<2q$ and $q\ge 2$. It is clearly that $b+q-1\le a\le b+q+1$. This implies $2a\le a+b+q+1\le 3q$ and $2b\le 2a-2q+2\le 3q-2q+2=q+2$. Hence $a\le \lfloor\frac{3q}{2} \rfloor$ and $b\le \min\{\lfloor\frac{q+2}{2} \rfloor, 2q-a-1\}$.


\ms\nt For $q=2$, since $a\ge 2$ and $b\ne 1$, we have $b=0$. So $(a,b)\in \{(2,0), (3, 0)\}$.  For $(a,b)=(2,0)$, the graph is not simple. For $(a,b)=(3,0)$, $G\cong P_3$. It is easy to verify that $k=1$ and the consecutive vertex labels must be $3,5,1$ or $3,1,5$.

\ms\nt Now we consider $q\ge 3$. For $b=0$, by Theorem~~\ref{thm-NAV}, the largest label must be odd. We have $q\le a\le q+1$. So $G$ is a $(q, q)$-graph or a $(q+1,q)$-graph.

\ms\nt For a $(q, q)$-graph with $q\ge 3$, we know that the largest label $k+2q-1$ is odd, i.e., $k$ is even. By Theorem~\ref{thm-NAV}, the second largest label $k+2q-2$ which is even must be a vertex label, contradicting $b=0$.

\ms\nt Examples of $(q+1,q)$-graphs with $q\ge 3$ that are caterpillar and spider graphs with all their vertices labeled by odd numbers can be found in~\cite[Theorems 4.1 \& 4.2]{Lau+Shiu+Ng}.

\begin{theorem}\label{thm-oddvertex} If $G$ is $k$-super graceful with all vertices labeled by odd integers, then $G$ is a $(q+1,q)$-graph and  $k=1$, where $q\ge 1$.\end{theorem}
\begin{proof}
Combining the above discussions, we know that $G$ is a $(q+1,q)$-graph, where $q\ge 1$. Since the number of odd labels is one more than that of even labels, $k$ is odd. Moreover, the maximum possible edge label is $(k+2q)-k=2q$.  Since $k+2q-1$ is an edge label, we have $k+2q-1\le 2q$ and hence $k=1$.
\end{proof}

\begin{theorem}\label{thm-SGnG} For a $(q+1,q)$-graph $G=(V,E)$, the following statements are equivalent:
\begin{enumerate}[{\rm (i)}]
  \item $G$ admits a super graceful labeling with edge-label set $\mb{[1,q]}$.
  \item $G$ is graceful.
  \item $G$ admits a super graceful labeling with all vertices labeled by odd integers.
\end{enumerate}
\end{theorem}

\begin{proof} (i) $\implies$ (ii). Let $f$ be a super graceful labeling of $G$ with edge-label set $\mb{[1,q]}$. For each $u\in V$, define $g(u)=f(u)-q-1$. Now, we can easily check that $g^*(uv)=f(uv)$ for each $uv\in E$. Since $f$ is bijection and $f(E)=\mb{[1,q]}$, $g^*$ is a bijection. Hence $g$ is a graceful labeling of $G$.

\ms\nt (ii) $\implies$ (iii). Let $g$ be a graceful labeling of $G$. For each $u\in V(G)$ and $uv\in E(G)$, define $h(u) = 2g(u)+1$ and $h(uv) = 2g^*(uv)$. Clearly, $h$ is a super graceful labeling of $G$ with all vertices labeled by odd integers.

\ms\nt (iii) $\implies$ (i). Let $h$ be a super graceful labeling of $G$ with odd vertex labels only. For each $u\in V(G)$ and $uv\in E(G)$, define $f(u) = (h(u)-1)/2 + q+1$ and $f(uv) = h(uv)/2$. Clearly, $f$ is a super graceful labeling of $G$ with edge-label set $\mb{[1,q]}$.
\end{proof}

\nt Let $G$ be a graph and $v\in V(G)$. Let $P_n=u_1\cdots u_n$ be a path of order $n\ge 2$. The graph $G^v_n$ is obtained from $G$ by merging $v$ with $u_1$. With this notation, we have the following.

\begin{theorem}\label{thm-extend} Suppose $G$ is a super graceful $(q+1,q)$-graph with all vertices labeled by odd integers, where $q\ge 1$. Let $v$ be the vertex labeled by $1$. There is a super graceful labeling of $G^v_n$ such that all vertices of this graph are labeled by odd integers and $u_n$ is labeled by $1$. \end{theorem}

\begin{proof}
Here we only show the case when $n=2$. For the general case, the proof is by mathematical induction on $n$.

\ms\nt Let $f$ be a super graceful labeling of $G$ such that all vertices of $G$ are labeled by odd integers and $f(v)=1$. Define a labeling $g$ on  $G^v_2$ by $g(w)=2q+4-f(w)$ for $w\in V(G)$, $g(u_2)=1$ and $g(xy)=|g(x)-g(y)|$. It is easy to check that $g$ is a  super graceful labeling.
\end{proof}

\begin{corollary}\label{cor-extend} If an $(n+1,n)$-graph $G$ is super graceful with all odd vertex labels, then $G$ admits two super graceful labelings $f$ and $g$ such that  $g(u) = 2n - f(u)$ for each $v\in V(G)$ and $g(uv) = f(uv)$ for each $uv\in E(G)$. \end{corollary}

\nt The following conjecture is equivalant to the Graceful Tree Conjecture.

\begin{conjecture} Each tree admits a super graceful labeling with all odd vertex labels.  \end{conjecture}

\begin{problem}\label{pbm-even} For $k\ge 1$, characterize $k$-super graceful graphs with all odd vertex labels. \end{problem}

\nt\section{\mathversion{bold}1- and  2-Regular graphs}

\nt We first consider $k$-super graceful $nK_2$ with all even edge labels.

\begin{lemma}\label{lem-nK2evenedge} If $nK_2$ is $k$-super graceful with all even edge labels, then $(i)$ $n\equiv3\pmod{4}$ and $k$ is even, or $(ii)$ $n\equiv0\pmod{4}$, or $(iii)$ $n\equiv1\pmod{4}$ and $k$ is odd.
\end{lemma}

\begin{proof} Assume that there are $s$ edges with odd end-vertex labels, where $s\ge 1$. Thus, there are $2(n-s)$ vertices with even labels. Hence, $-1\le 2s - (3n-2s)\le 1$. If $4s- 3n = -1$, then $n\equiv3\pmod{4}$ and $k$ is even since we need an odd number of even labels. If $4s - 3n = 0$, then $n\equiv0\pmod{4}$. If $4s - 3n=1$, then $n\equiv1\pmod{4}$ and $k$ is odd since we need an odd number of even labels.
\end{proof}

\nt It follows that (also by Corollary~\ref{cor-1SG}), $nK_2$ is super graceful with all even edge labels only if $n\equiv0,1\pmod{4}$. Suppose such a labeling for $nK_2$ exists, we must have triples $(a_i,b_i,c_i)$ with $b_i = c_i - a_i$ such that the $b_i$'s are distinct even numbers and the $a_i$'s and $c_i$'s are distinct odd numbers from 1 to $3n$. If $n=4t\ge 4$, then $1\le i\le 3s$ and we are left with $3s$ even numbers, say $d_j, e_j, f_j$ ($1\le j\le s$) such have $f_j = d_j + e_j$. However, if $\sum f_j \not= \sum (d_j + e_j)$, then such a labeling does not exist. Similarly, we can also find such a labeling for $n=4t+1\ge 5$. Observe that $\sum (d_j+e_j+f_j) \equiv0\pmod{4}$.

\ms\nt For $n=4,5,8,9$, the following labelings show that $nK_2$ is super graceful with all their edges labeled by even numbers.
\begin{enumerate}[(1)]
  \item $n=4$: \begin{enumerate}[a.]
  \item Edge labels are $2,4,6,10$; adjacent vertex labels are $(5,7)$, $(8,12)$, $(3,9)$, $(1,11)$. So, the corresponding triples $(a_i, b_i, c_i)$'s are $(5,2,7)$, $(3,6,9)$, $(1,10,11)$ and $(d_1, e_1, f_1)$ is $(8,4,12)$.
      \item Edge labels are $2,4,6,8$; adjacent vertex labels are $(10,12)$, $(5,9)$, $(1,7)$, $(3,11)$. So, the corresponding triples $(a_i, b_i, c_i)$'s are $(5,4,9)$, $(1,6,7)$, $(3,8,11)$ and $(d_1, e_1, f_1)$ is $(10,2,12)$.
      \end{enumerate}
  \item $n=5$: Edge labels are $2,4,6,10,14$; adjacent vertex labels are $(7,9)$, $(8,12)$, $(5,11)$, $(3,13)$, $(1,15)$. So, the corresponding triples $(a_i, b_i, c_i)$'s are $(7,2,9)$, $(5,6,11)$, $(3,10,13)$, $(1,14,15)$ and $(d_1, e_1, f_1)$ is $(8,4,12)$.
  \item $n=8$: Edge labels are $2,4,6,8,10,12,14,18$; adjacent vertex labels are $(20,22)$, $(7,11)$, $(15,21)$, $(16,24)$, $(3,13)$, $(5,17)$, $(9,23)$, $(1,19)$. So, the corresponding triples $(a_i,b_i,c_i)$'s are $(7,4,11)$, $(15,6,21)$, $(3,10,13)$, $(5,12,17)$, $(9,14,23)$, $(1,18,19)$ and $(d_i,e_i,f_i)$'s are $(20,2,22)$, $(16,8,24)$.
  \item $n=9$: Edge labels are $2,4,6,8,10,12,18,20,24$; adjacent vertex labels are $(14,16)$, $(22,26)$, $(9,15)$, $(11,19)$, $(7,17)$, $(13,25)$, $(5,23)$, $(1,21)$, $(3,27)$. So, the corresponding triples $(a_i,b_i,c_i)$'s are $(9,6,15)$, $(11,8,19)$, $(7,10,17)$, $(13,12,25)$, $(5,18,23)$, $(1,20,21)$, $(3,24,27)$ and $(d_i,e_i,f_i)$'s are $(14,2,16)$, $(22,4,26)$.
\end{enumerate}

 \begin{theorem} For $t\ge 1$, if $(4t)K_2$ is super graceful with all even edge labels, then so is $(4t+1)K_2$. \end{theorem}

\begin{proof} Write  $(4t)K_2=G+H$, where $G\cong (3t)K_2$ and $H\cong tK_2$ such that the vertices of $G$ are labeled by odd numbers and the vertices of $H$ are labeled by even numbers. Now we label the vertices of the graph $(4t+1)K_2\cong G+H+K_2$ as follows:

\ms\nt Add 2 to each of the original labels of the vertices in $G$; keep the labels of $H$; and label the vertices of the last $K_2$ by $1$ and $12t+3$.

\ms\nt Now, odds from $3$ to $12t+1$ appear in $G$; evens from $2$ to $12t$ appear in $G+H$. Only $1$, $12t+2$ and $12t+3$ are not involved in $G+H$, but they appear in $K_2$. Hence $(4t+1)K_2$ is super graceful. \end{proof}

\nt It is clear that $3K_2$ is not 2-super graceful with all even edge labels. Hence, the conditions of Lemma~\ref{lem-nK2evenedge} may not be sufficient.

\begin{problem} Determine which necessary condition in Lemma~\ref{lem-nK2evenedge} is also sufficient. \end{problem}

\begin{lemma}~\label{lem-kSGkSA} For $n\ge 1$, $nK_2$ is $k$-super graceful if and only if it is $k$-sequentially additive.    \end{lemma}

\begin{proof} Let $f$ be a $k$-super graceful labeling of $nK_2$. For each $K_2$ component with vertex labels $f(u) > f(v)$ and edge label $f(uv)$, we have $f(u)-f(v)=f(uv)$ if and only if $f(v) + f(uv) = f(u)$. Define $g$ such that $g(v)=f(v)$, $g(u) = f(uv)$ and $g(uv) = f(u)$. Clearly $g$ is a $k$-sequentially additive labeling of $nK_2$.  By a similar argument, we have the converse.  \end{proof}

\nt From the proof of Theorem 3.1 in~\cite{Hegde+Miller}, we have

\begin{lemma}\label{lem-nonkSA} (i) If $n\equiv 2$ or $3\pmod{4}$, then $nK_2$ is not $k$-sequentially additive for odd $k\ge 1$. (ii) If $n\equiv 2\pmod{4}$, then $nK_2$ is not $k$-sequentially additive for even $k\ge 2$.  \end{lemma}

\nt Lemma~\ref{lem-nonkSA} implies that $nK_2$ is not $n$-sequentially additive for $n\equiv2,3\pmod{4}$. In fact,

\begin{theorem}\label{thm-nSG} The graph $nK_2$ is $n$-super graceful with $n\ge 1$, and hence $n$-sequentially additive, if and only if  $n=1$.  \end{theorem}

\begin{proof} By Corollary~\ref{cor-NAV}, we must label $nK_2$ such that for each of the $n$ copies of $K_2$, an end-vertex is labeled by an integer in $\mb{[3n,4n-1]}$ with the property that these $n$ labels are distinct. Hence, we must form $n$ pairs of integers in $\mb{[n,3n-1]}$ such that the sum of the pairs are distinct integers from $3n$ to $4n-1$. Therefore, we must have $\sum^{3n-1}_{i=n} i = \sum^{4n-1}_{i=3n} i$ if and only if $n=1$.
\end{proof}

\begin{definition} For $n,k\ge 1$, a {\it $k$-Skolem sequence of length $n$} is a sequence $\{(a_i,b_i)\}^n_{i=1}$ such that $a_i-b_i=k+i-1$ for mutually distinct $a_i,b_i\in \mb{[1,2n]}$.  \end{definition}

\nt A 1-Skolem sequence is also called a Skolem sequence (see~\cite{Skolem}).

\begin{lemma}[\cite{Skolem}]\label{lem-skolem} A Skolem sequence of length $n$ exists if and only if $n\equiv 0,1\pmod{4}$.  \end{lemma}

\begin{corollary}\label{cor-1SG} $(i)$ For $n\ge 1$, $nK_2$ is super graceful if and only if $n\equiv 0,1\pmod{4}$. Moreover, it is possible to have $\mb{[1,n]}$ as the edge-label set.  $(ii)$ For $n\ge 2$, $nK_2$ is not $k$-super graceful for odd $k\ge 1$ if $n\equiv 2$ or $3\pmod{4}$, and for even $k\ge 2$ if $n\equiv 2\pmod{4}$. Consequently, $nK_2$ is not $k$-super graceful for all $k$ if $n\equiv 2\pmod{4}$. $(iii)$ For all $n\ge 2$, $nK_2$ is not $n$-super graceful.   \end{corollary}

\begin{proof} (i) It follows from Lemmas~\ref{lem-kSGkSA}, \ref{lem-nonkSA} and~\ref{lem-skolem} by extending a Skolem sequence $(a_i,b_i)$ to sequence of triples $(a_i+n,b_i+n,i)$ ($1\le i\le n$) that induces a super graceful labeling of $nK_2$ with edge-label set $\mb{[1,n]}$. (ii) and (iii) follow from Lemma~\ref{lem-nonkSA} and Theorem~\ref{thm-nSG}, respectively.   \end{proof}




\begin{lemma}\label{lem-S-SG} A $k$-Skolem sequence of length $n$ exists if and only if $nK_2$ is $k$-super graceful  with edge-label set $\mb{[k,k+n-1]}$.
\end{lemma}

\begin{proof} (Necessity) Let $\{(a_i,b_i)\}^n_{i=1}$ be a $k$-Skolem sequence. By adding $k+n-1$ to each $a_i$ and $b_i$, we get a $k$-super graceful labeling of $nK_2$  with edge-label set $\mb{[k,k+n-1]}$.

\ms\nt (Sufficiency) A $k$-super graceful labeling of $nK_2$  with edge-label set $\mb{[k,k+n-1]}$ has vertex-label set $\mb{[k+n,k+3n-1]}$. Subtracting $k+n-1$ from each vertex label gives a $k$-Skolem sequence of length $n$. \end{proof}


\begin{theorem}\label{thm-kSkol} If a $k$-Skolem sequence of length $n$ exists, then $n\equiv0,3\pmod{4}$ for even $k$, and $n\equiv0,1\pmod{4}$ for odd $k$. \end{theorem}

\begin{proof} Suppope a $k$-Skolem sequence of length $n$ exists, then $a_i-b_i = k+i-1$ for $1\le i\le n$. Thus, we have \[\sum a_i - \sum b_i = \frac{n(n+2k-1)}{2}.\] On the other hand, \[\sum a_i + \sum b_i = n(2n+1).\] Hence, adding the two equations yields \[\sum a_i = \frac{n(5n+2k+1)}{4}\] which is an integer whenever $n\equiv0,3\pmod{4}$ for even $k$, and $n\equiv0,1\pmod{4}$ for odd $k$.
\end{proof}

\nt For $1\le i\le n$, let $e_i = u_iv_i$ be the edge of the $i$-th component of $nK_2$.

\begin{lemma}\label{lem-nK2NC} Suppose $k,n\ge 1$. If $nK_2$ is $k$-super graceful, then $3n(2k+3n-1)/2$ is even and $n\ge 2k-1$.
\end{lemma}

\begin{proof} Since $nK_2$ is $k$-super graceful, the labels are from $k$ to $k+3n-1$. Without loss of generality, assume that $f(u_i) = f(v_i) + f(e_i)$ for $1\le i\le n$. Hence, $\sum_{i=1}^n f(u_i) = \sum_{i=1}^n f(v_i) + \sum_{i=1}^n f(e_i)$. Therefore, $\sum_{i=1}^n f(u_i) + \sum_{i=1}^n f(v_i) + \sum_{i=1}^n f(e_i) = 3n(2k+3n-1)/2$ is even. This further implies that the sum of the largest $n$ labels is greater than or equal to $3n(2k+3n-1)/4$. Consequently, $\sum_{i=1}^n 2n+k-1+i = n(5n+2k-1)/2 \ge 3n(3n+2k-1)/4$ giving us $n\ge 2k-1$.  \end{proof}


\nt Applying the contrapositive of Lemma~\ref{lem-nK2NC} we have
\begin{corollary} The graph $nK_2$ is not $k$-super graceful for (i) $n\equiv 2\pmod{4}$ and all $k$; (ii) $n\equiv1\pmod{4}$ and even $k$; (iii) $n\equiv3\pmod{4}$ and odd $k$. \end{corollary}

\begin{theorem}\label{thm-kSGnK2} For $n,k\ge 1$, $nK_2$ is $k$-super graceful with the largest $n$ labels assigned to mutually non-adjacent vertices if and only if $n=2k-1$. Moreover, $(2k-1)K_2$ is $k$-super graceful with edge-label set $\mb{[k,3k-2]}$.\end{theorem}

\begin{proof} (Sufficiency) Suppose $nK_2$ is $k$-super graceful with the largest $n$ labels assigned to mutually non-adjacent vertices. We then have $\sum_{i=k}^{k+2n-1} i = \sum_{j=k+2n}^{k+3n-1} j$. Solving the equation gives $n=2k-1$.

\ms\nt (Necessity) It suffices to find a $k$-Skolem sequence $\{(a_i, b_i)\}$ of length $2k-1$. It is presented below:
$$\begin{array}{l|c|c|l}
\mbox{range of }i  & a_i & b_i & a_i-b_i\\\hline
1\le i \le k & 2k-2+2i & k-1+i & \multirow{2}{1.5cm}{$k-1+i$}\\
k+1\le i \le 2k-1 & 2i-1 & i-k &
\end{array}$$
\nt By Lemma~\ref{lem-S-SG}, $(2k-1)K_2$ is $k$-super graceful with edge-label set $\mb{[k,3k-2]}$. 
\end{proof}

 \begin{corollary}\label{cor-special} If $nK_2$ is $k$-super graceful with edge-label set $\mb{[k, k+n-1]}$, then $(3n+2k-1)K_2$ is also $k$-super graceful with edge-label set $\mb{[k, 3n+3k-2]}$.\end{corollary}

\begin{proof}
By Theorem~\ref{thm-kSGnK2}, we know that $(2n+2k-1)K_2$ is $(n+k)$-super graceful with edge-label set $\mb{[n+k, 3n+3k-2]}$. By Lemma~\ref{lem-G+H4} $(3n+2k-1)K_2$ is $k$-super graceful with edge-label set $\mb{[k, 3n+3k-2]}$.
\end{proof}

\begin{corollary}\label{cor-2SGnK2} If $n=(2k-1)(3^r-1)/2$ for $r\ge 1$, then $nK_2$ is $k$-super graceful with edge-label set $\mb{[k,k+n-1]}$.  \end{corollary}

\begin{proof} The corollary follows by applying Corollary~\ref{cor-special} and mathematical induction on $r$.
\end{proof}

\begin{corollary}\label{cor-kSS} (i) For each $k\ge 1$, there exist a $k$-Skolem sequence of length $2k-1$. (ii) For each $k\ge 1$, $(2k-1)K_2$ is $k$-sequentially additive with edge-label set $\mb{[5k-2,7k-4]}$.  (iii) For each $k\ge 1$, there is a $k$-super graceful graph $G$ of size $n = 2k-1$ with edge-label set $\mb{[k,k+n-1]}$. Particularly, for each $k\ge 1$, there are infinitely many $k$-super graceful $1$-regular graphs. \end{corollary}

\begin{theorem}\label{thm-2SGnK2}  Suppose $n\le 84$. $nK_2$ is $2$-super graceful with edge-label set $\mb{[2,n+1]}$ if and only if $n\equiv0,3\pmod{4}$. \end{theorem}

\begin{proof} The necessity follows from Lemma~\ref{lem-S-SG} and Theorem~\ref{thm-kSkol}. We prove the sufficiency by constructing a $2$-Skolem sequence of length $n$.

\ms\nt Let $n=4r-1$ for $1\le r\le 21$.\\
For $r=1$, a $2$-Skolem sequence of length $3$ is given by
$(4,2)$, $(6,3)$ and $(5,1)$.\\
For $r=2$, a $2$-Skolem sequence of length $7$ is given by
$(10,8)$, $(4,1)$, $(6,2)$, $(14, 9)$, $(13, 7)$, $(12, 5)$, $(11,3)$.\\
For $r=3$, a $2$-Skolem sequence of length $11$ is given by $(14,12)$, $(4,1)$, $(6,2)$, $(8,3)$, $(16, 10)$, $(22, 15)$, $(21, 13)$, $(20, 11)$, $(19,9)$,
$(18,7)$, $(17,5)$.\\
$$(a_i, b_i)=\begin{cases}
(2i, i-1) & \mbox{for }2\le i\le 4 \mbox{ or even $i$ in the range } 6\le i \le 2r-2;\\
(10r-2-i, 10r-3-2i) &  2r\le i \le 4r-1.
\end{cases}$$
It is easy to see that the difference $a_i-b_i=i+1$.
Up to now, the values of $a_i$'s cover the range $\{4,6,8\}\cup\{4s\;|\; 3\le s\le r-1\}\cup \mb{[6r-1, 8r-2]}$, the values of $b_i$'s cover the range $\{1,2,3\}\cup\{2s-1\;|\; 3\le s\le 3r-1\}$. That is, the set of unassigned integers is $\{4l+2\;|\; 2\le l \le r-2\}\cup\{2l\;|\; 2r-1\le l\le 3r-1\}$ which will be assigned to $a_i$ and $b_i$ for odd integers $i$ in the range $\{1\}\cup \mb{[5, 2r-1]}$.
The difference $(a_i-b_i)$'s cover the range $\{3,4,5\}\cup \{2s-1\;|\; 4\le s\le r\}\cup \mb{[2r+1, 4r]}$. That is, the set of missing differences is $\{2\}\cup \{2l\;|\; 3\le l\le r\}$. Ad hoc assignments for the remaining $(a_i, b_i)$'s are shown below.
$$\begin{array}{l|l|l}
r & (a_1, b_1) & (a_i, b_i) \mbox{ for odd $i$ in \mb{[5, 2r-1]}}\\\hline
4 & (20,18) & (16,10), (22,14)\\
5 & (28,26) & (24,18), (22,14), (20,10)\\
6 & (32,30) & (34,28), (26,18), (24,14), (22,10)\\
7 & (36,34) & (38,32), (26,18), (40,30), (22,10), (28,14)\\
8 & (42,40) & (38,32), (18,10), (44,34), (26,14), (36,22), (46,30)\\
9 & (52,50) & (48,42), (46,38), (44,34), (22,10), (40,26), (30,14), (36,18)\\
10 & (42,40) & (44,48), (18,10), (56,46), (26,14), (58,44), (38,22), (52,34), (50,30)\\
11 & (64,62) & (54,48), (58,50), (52,42), (30,18), (60,46), (26,10), (56,38), (34,14), (44,22)
\end{array}$$

\nt Using computer searching, we have also obtained at least one such pairings for $12\le r\le 21$.  In fact, our computer search shows that there are altogether 189 different pairings for $r=11$. For $r=20$, there are at least 5657 different pairings. 

\ms\nt Let $n=4r$ for $1\le r\le 21$. From the $2$-Skolem sequence $\{(a_i, b_i)\}_{i=1}^{4r-1}$ of length $4r-1$ defined above, we define
$$(A_i, B_i)=\begin{cases}
(a_i, b_i) & \mbox{if } 1\le i\le 2r-1,\\
(a_i+2, b_i+2) & \mbox{if } 2r\le i\le 4r-1,\\
(6r, 2r-1) & \mbox{if } i=4r.
\end{cases}$$
One may easily check that $\{(A_i, B_i)\}_{i=1}^{4r}$ is a $2$-Skolem sequence of length $4r$.
%
\end{proof}

\nt  Clearly, $nK_2$ is 2-super graceful for all $n\equiv0,3\pmod{4}$ if the following conjecture holds. Moreover, for $n\equiv 0\pmod{4}$, we can get a labeling with all even edge labels.

\begin{conjecture} For $r\ge 22$, there exist $r-1$ pairs of integers using $\{4l+2\;|\; 2\le l \le r-2\}\cup\{2l\;|\; 2r-1\le l\le 3r-1\}$ completely such that the set of differences of the pairs is $\{2\}\cup \{2l\;|\; 3\le l\le r\}$.\end{conjecture}

\nt By  Corollary~\ref{cor-special} and Theorem~\ref{thm-2SGnK2}, we obtain infinitely many 2-super graceful 1-regular graphs recursively.  By Lemma~\ref{lem-nK2NC}, it is easy to show that $4K_2$ is not 3-super graceful but $8K_2$ is 3-super graceful given by the triples  $(10,3,13)$, $(18,4,22)$,  $(16,5,21)$, $(14,6,20)$, $(19,7,26)$, $(17,8,25)$, $(15,9,24)$, $(12,11,23)$. Consequently, $3K_2$ is $k$-super graceful if and only if $k=2$; $4K_2$ is $k$-super graceful if and only if $k=1,2$; and $5K_2$ is $k$-super graceful if and only if $k=1,3$.

\ms\nt We can also show that $nK_2$ is 3-super graceful for $n=8,9,12,13$. For $n=8$, the triples are $(13,3,16)$, $(14,4,18)$, $(21,5,26)$, $(19,6,25)$, $(17,7,24)$, $(15,8,23)$, $(11,9,20)$, $(12,10,22)$. For $9K_2$, the triples are $(24,5,29)$, $(22,6,28)$, $(20,7,27)$, $(17,9,26)$, $(15,10,25)$, $(12,11,23)$, $(13,8,21)$, $(16,3,19)$, $(14,4,18)$.  For $n=12$, the triples are  $(18,3,21)$, $(19,4,23)$, $(20,5,25)$, $(32,6,38)$, $(30,7,37)$, $(28,8,36)$, $(26,9,35)$, $(24,10,34)$, $(22,11,33)$, $(15,12,27)$, $(16,13,29)$ , $(17,14,31)$  (where the last three triples may be replaced by $(14,13,27)$, $(17,12,29)$, $(16,15,31)$). For $13K_2$, the triples are $(36,5,41)$, $(34,6,40)$, $(32,7,39)$, $(30,8,38)$, $(22,15,37)$, $(21,14,35)$, $(20,13,33)$, $(19,12,31)$, $(18,11,29)$, $(24,4,28)$, $(17,10,27)$, $(23,3,26)$, $(16,9,25)$. Note that each labeling above gives the edge-label set $\mb{[3, n+2]}$. Moreover, $12K_2$ is 5-super graceful with triples $(32,8,40)$, $(30,9,39)$, $(28,10,38)$, $(22,15,37)$, $(20,16,36)$, $(18,17,35)$, $(23,11,34)$, $(21,12,33)$, $(25,6,31)$, $(24,5,29)$, $(14,13,27)$, $(19,7,26)$ but the edge-label set is not $\mb{[5,16]}$.

\begin{conjecture} Every $k$-super graceful $nK_2$ admits a labeling with edge-label set $\mb{[k,k+n-1]}$.  \end{conjecture}

\ms\nt In the rest of this section, we discuss 2-regular graphs. We note that $C_{10}$ with 5-super graceful with consecutive vertex labels 24, 6, 20, 12, 21, 10, 23, 7, 22, 5 and corresponding edge labels 18, 14, 8, 9, 11, 13, 16, 15, 17, 19. Moreover, $C_{12}$ is 6-super graceful with consecutive vertex labels 27, 6, 29, 7, 26, 14, 24, 13, 28, 8, 25, 9 and corresponding edge labels 21, 23, 22, 19, 12, 10, 11, 15, 20, 17, 16, 18.

\begin{problem} For $k\ge5$, find sufficient conditions such that $C_{2k}$ is $k$-super graceful. \end{problem}









\nt The following follows directly from the definition.

\begin{lemma}~\label{lem-edgeC3} In any $k$-super graceful labeling of a graph, the sum of the two smaller edge labels of a $C_3$ subgraph is the label of the third edge.   \end{lemma}

\nt  In~\cite{Hajnal+Nagy}, the authors showed that 2-regular graphs $nC_3$ is SSA for all $n\ge 1$ and $nC_4$ is SSA if and only if $n\equiv 0,1\pmod{3}$. However,

\begin{theorem}~\label{thm-kC3}  For $k\ge 1$, $kC_3$ is $k$-super graceful if and only if $k=1$.      \end{theorem}

\begin{proof} Let $u_{i,1}, u_{i,2}$ and $u_{i,3}$ be the three vertices of the $i$-th component of $kC_3$. Since every cycle is super graceful, the sufficiency holds. To prove the necessity, we may assume $k\ge 2$. Suppose $kC_3$ admits a $k$-super graceful labeling $f$. By Theorem~\ref{thm-NAV}, without loss of generality, we may assume that $f(u_{i,1})=7k-i$, $1\le i\le k$. Under this condition, $6k-1$ must a vertex or an edge label in the first component, i.e., (i) $f(u_{1,2})=6k-1$ or (ii) $f(u_{1,1}u_{1,2})=6k-1$.

\ms\nt If (i) holds, then $f(u_{1,1}u_{1,2})=k$. It follows that $f(u_{1,1}u_{1,3})=6k-2$. By Lemma~\ref{lem-edgeC3} $f(u_{1,2}u_{1,3})=5k-2$ and hence $f(u_{1,3})=k+1$. Now $6k-3$ is neither a vertex nor an edge label, a contradiction.

\ms\nt If (ii) holds, then $f(u_{1,2})=k$. It follows that $f(u_{1,3})=6k-2$ and $f(u_{1,1}u_{1,3})=k+1$. Consequently, $f(u_{1,2}u_{1,3})=5k-2$ and $6k-3$ is neither a vertex nor an edge label, a contradiction.
\end{proof}

\begin{corollary}  The graph $2C_3$ is not $k$-super graceful for each $k\ge 1$. \end{corollary}

\begin{proof} By Theorem~\ref{thm-kC3}, $2C_3$ is not $2$-super graceful. By Lemma~\ref{lem-edgeC3}, we can show that $2C_3$ is not super graceful. \end{proof}

\section{Bi-regular graphs}

\nt It is known that $P_n, n\ge 2$, is super graceful with all even edge labels.  Moreover, $nK_2$ is super graceful with all even edge labels for $n=4,5,8,9$.

\begin{problem} For which $k\ge 2$ is there a $k$-super graceful graph with all even edge labels? \end{problem} 

\begin{theorem}~\label{thm-oddedge} A graph $G$ is $k$-super graceful with all odd edge labels if and only if $k=1$ and $G$ is a star $K(1,q)$, $q\ge 1$.
\end{theorem}

\begin{proof} (Sufficiency) Suppose $G=K(1,q)$. We first label the central vertex of $G$ by $2q+1$. We then label the edges by $1,3,5,\ldots, 2q-1$, the corresponding end-vertices by $2q, 2q-2, 2q-4, \ldots, 2$, respectively. We have a required labeling.

\ms\nt (Necessity) Observe that any two adjacent vertex labels must have different parity. Hence, $G$ is bipartite. Suppose $G$ has $q$ edges. Let $G$ have partite sets $A,B$ such that all vertices in $A$ have odd labels and all vertices in $B$ have even labels.  Since $q\ge 1$, $|A|=a\ge 1$ and $|B|\ge 1$.  Suppose $k$ is odd. We have $|B|=a+q-i\ge 1$, $i=0,1$. Note that $|B|\le q$ which implies that $a+q-i\le q$. So $a\le i$. Hence, $a=i=1$. Consequently, $G$ must be a star $K(1,q), q\ge 1$ with largest label $2q+k$ which is odd. If $k\ge 3$, by Theorem 4.2 in~\cite{Lau+Shiu+Ng-kSG}, we know $k$ is the central vertex label of $K(1,q)$.  This means $2q+k$ must be an edge label, contradicting Theorem~\ref{thm-NAV}. Hence, $k=1$.

\ms\nt Suppose $k$ is even. Now $|B|=a+q+i\ge 1$, $i=0,1$. Again, $|B|\le q$ implies that $a+i\le 0$, a contradiction. \end{proof}

\nt It was shown in~\cite{Perumal1, Lau+Shiu+Ng-kSG} the complete bipartite graph $K(m,n)$ is $1$-, $m$- and $n$-super graceful for all $n\ge m\ge 2$ whereas $K(1,q)$ is $k$-super graceful if and only if $q\equiv0\pmod{k}$.

\begin{corollary} For each $k\ge 1$, there exists a $k$-super graceful $G$ such that $\Delta(G)=k$.  \end{corollary}

\nt We now study the $k$-super gracefulness of the bi-regular graph $K(m,n)$. Let $A,B$ be the two partite sets of $K(m,n)$ with $|A|=m, |B|=n$. Without loss of generality, we may assume that a $k$-super graceful labeling $f$ of $K(m,n)$ assigns even integers to $a$ vertices in $A$ and $b$ vertices in $B$. Hence, the total number of even labels under $f$ is $a+b+ab + (m-a)(n-b)$, and the total number of odd labels under $f$ is $(m-a) + (n-b) + b(m-a) + a(n-b)$.

\ms\nt Suppose $m,n$ are even. We then have equal numbers of odd and even labels. Hence, $(m-a) + (n-b) + b(m-a) + a(n-b) = a+b+ab + (m-a)(n-b)$ which implies that $(m-2a-1)(n-2b-1)=1$ so that $a=m/2, b=n/2$ or $a=(m-2)/2, b=(n-2)/2$.

\ms\nt Suppose $m$ or $n$ is odd. We then consider 2 cases:

\ms\nt Case (1): $k$ is odd. In this case, the number of even labels is one less than the number of odd labels. Hence, $a+b+ab + (m-a)(n-b) + 1 = (m-a) + (n-b) + b(m-a) + a(n-b)$ which implies that $(m-2a-1)(n-2b-1)=0$ so that $a=(m-1)/2$ or $b=(n-1)/2$.

\ms\nt Case (2): $k$ is even. In this case, the number of odd labels is one less than the number of even labels. Hence, $a+b+ab + (m-a)(n-b) = (m-a) + (n-b) + b(m-a) + a(n-b) + 1$ which implies that $(m-2a-1)(n-2b-1)=2$ so that $a=(m-2)/2, b=(n-3)/2$ or $a=(m-3)/2, b=(n-2)/2$ or $a=m/2, b=(n+1)/2$ or $a=(m+1)/2, b=n/2$. Consequently, $m$ and $n$ must be of different parity.

\begin{theorem}\label{thm-Kmn} Suppose $K(m,n)$ is $k$ super graceful, where $n\ge m\ge 2$, and that $a$ and $b$ are as defined above.
\begin{enumerate}[(i)]
  \item If $m,n$ are even, then $a=m/2, b=n/2$ or $a=(m-2)/2, b=(n-2)/2$.
  \item If $m,n$ are odd, then $k$ is odd with $a=(m-1)/2$, $b=(n-1)/2$.
  \item If $m,n$ are of different parity, then $k$ is odd implies that $m$ is odd with $a=(m-1)/2$ or $n$ is odd with $b=(n-1)/2$; and $k$ is even implies $a=(m-2)/2, b=(n-3)/2$; or $a=(m-3)/2, b=(n-2)/2$; or $a=m/2, b=(n+1)/2$; or $a=(m+1)/2, b=n/2$.
\end{enumerate}
\end{theorem}

\nt By Theorem~\ref{thm-Kmn}, we conclude that $K(3,3)$ is not 2-super graceful, $K(2,4)$ is not 3-super graceful and $K(3,5)$ is  neither 2- nor 4-super graceful. Suppose $K(3,4)$ has a 2-super graceful labeling $f$. By Theorem~\ref{thm-Kmn}~(iii), we have $a=0, b=1$ or $a=b=2$.  Suppose $a=0, b=1$. By Theorem~\ref{thm-NAV}, 20 and 19 are labels of vertices of $B$. Hence, 18 cannot be a label, a contradiction. Suppose $a=b=2$. Let $A=\{x,y,z\}$ and $B=\{s,t,u,v\}$.  Suppose $f(x), f(y), f(s), f(t)$ are even. By Theorem~\ref{thm-NAV}, we consider two cases.

\begin{enumerate}[(a)]
\item $f(x)=20,f(z)=19$. We may assume (1) $f(xs)=18$ or (2) $f(y) = 18$.

In (1), we have $f(s)=2$ and $f(zs)=17$.  This forces, without loss of generality, (1.1) $f(yu)=3$ or (1.2) $f(zt)=3$.
\begin{enumerate}[{(1.}1)]
\item $f(yu)=3$ implies $f(xt)=16$ or $f(y)=16$. When $f(xt)=16$ we have $f(t)=4$ and $f(zt)=15$. This forces $f(yv)=5$ and $f(y)=14$. Hence $f(v)=9$. Now, $f(zv)=10=f(yt)$, a contradiction.  When $f(y)=16$, $f(ys)=14$, $f(u)=13$, $f(zu)=6$, and $f(xu)=7$. This forces $f(zv)=10$, $f(v)=9$, so that $f(yv)=7$, a contradiction.
\item $f(zt)=3$ implies $f(t)=16$ and hence $f(xt)=4$. This forces, without loss of generality, $f(zu)=6$. So $f(u)=13$ and $f(xu)=7$. Now, $15$ must be assigned to $xv$. Hence $f(v)=5$ and $f(zv)=14$. Now, there is no way to assign 12.
\end{enumerate}

In (2), without loss of generality, we may assume $f(ys)=2$ so that $f(s)=16$, $f(xs)=4$ and $f(zs)=3$. Now, there is noway to assign 17.

\item $f(s)=20,f(u)=19$. We may assume (1) $f(xs)=18$ or (2) $f(t)=18$. In (1), we have $f(x)=2$ and $f(xu)=17$. This forces (1.1) $f(t)=16$ or (1.2) $f(ys)=16$.
\begin{enumerate}[{(1.}1)]
\item $f(t)=16$ implies $f(xt)=14$. Then there is no room for the label 15.
\item $f(ys)=16$ implies $f(y)=4$ and $f(yu)=15$. This forces $f(t)=14$. Hence $f(xt)=12$ and $f(yt)=10$. It is easy to check that $f(zu)$ is neither 6 nor 8. So $f$ does not exist.
\end{enumerate}


In (2), we have (2.1) $f(zs)=17$ or (2.2) $f(v)=17$.
\begin{enumerate}[{(2.}1)]
\item $f(zs)=17$ implies $f(z)=3$, $f(zu)=16$ and $f(zt)=15$. This forces that $f(x)=14$ or, without loss of generality, $f(xs)=14$.     When $f(x)=14$, there is no room to assign the label 13.    When $f(xs)=14$, we have $f(x)=6$, $f(xt)=12$ and $f(xu)=13$. Since $(f(yt), f(yu), f(ys))$ is a sequence of three consecutive undetermined labels starting with an even label, $(f(yt), f(yu), f(ys))=(8,9,10)$. But this forces $f(y)=10$ which is impossible.

\item $f(v)=17$ implies, without loss of generality,  $f(xs)=16$. So that $f(x)=4$, $f(xt)=14$, $f(xu)=15$, $f(xv)=13$. Now $(f(yv), f(yt), f(yu), f(ys))$ and $(f(zv), f(zt), f(zu), f(zs))$ are two disjoint sequence of four consecutive undetermined labels. The first one starts with an odd label and the last one starts with an even label. But they  do not exist.
\end{enumerate}
%
%
%
%
%
\end{enumerate}

\ms\nt Hence, for $(m,n)= (2,2)$, $(2,3)$, $(2,4)$, $(3,3)$, $(3,4)$, $(3,5)$, $K(m,n)$ is $k$-super graceful if and only if $k=1,m,n$.

\begin{lemma} Suppose $n\ge 2$, the graph $K(2,n)$ is $k$-super graceful if and only if $k=1,2,n$. \end{lemma}

\begin{proof} Let the two partite sets of $K(2,n)$ be $A = \{u_1, \ldots, u_n\}$ and $B = \{v_1, v_2\}$.  Furthermore, let $E_1$ be the set of edges joining $v_1$ to the vertices in $A$, and $E_2$ be the set of edges joining $v_2$ to the vertices in $A$. The sufficiency follows from~\cite[Theorem 2.5]{Perumal1} and~\cite[Theorem 4.5]{Lau+Shiu+Ng-kSG}. To prove the necessity, we just need to show that if $K(2, n)$ is $k$-super graceful, with  $k \ge 3$ and $n\ge 3$, then $k = n$.

\ms\nt Let $f$ be a $k$-super graceful labeling of $K(2,n)$. The available label set is $\mb{[k, k + 3n + 1]}$. We may assume $f(u_1)>f(u_2)>\cdots >f(u_n)$ and $f(v_1)<f(v_2)$.  Let $j$ be the greatest so that the $j$ integers $k + 3n + 1, k + 3n, \ldots, k + 3n + 2 - j$ are labels of mutually non-adjacent vertices.  By Theorem~\ref{thm-NAV}, we have $n \ge j \ge k \ge 3$.  Thus these vertices are in $A$. Hence $f(u_i) = k + 3n + 2 - i$, for $i = 1, \ldots, j$.  Consider the largest undetermined label, $k + 3n + 1 - j$.  It cannot label a vertex in $A$, by the maximality of $j$. If it labels a vertex in $B$, then the edge joining it to $u_j$ has label $1 < k$, which is impossible.  Thus it must be the largest edge label. So we must have $f(u_1v_1) = k + 3n + 1 - j$, $f(v_1)=j$. Thus $f(u_iv_1) = k + 3n + 2 - i - j$, for $i = 1, \ldots, j$.
Consider the label $k + 3n + 1 - 2j$. Since it is the largest undetermined label, it must be the label of $u_1v_2$, $v_2$ or $u_{j+1}$ if it exists.
%

\ms\nt Suppose $f(u_1v_2) = k + 3n + 1 - 2j$. We have $f(v_2)= 2j$. Consequently, labels of $u_2v_2$ to $u_jv_2$ are $k+3n-2j$ to $k+3n+2-3j$.  If $u_{j+1}$ exists, then we must label $u_{j+1}$ to $u_{2j}$ by $k+3n+1-3j$ to $k+3n+2-4j$ so that the labels of $u_{j+1}v_2$ to $u_{2j}v_2$ are $k+3n+1-5j$ to $k+3n+2-6j$. This argument can be repeated until all vertices in $A$ are labeled. Hence $n\equiv 0\pmod{j}$. We see that there is a gap of integers from $j+1$ to $2j-1$ and all other used labels, except $j$ and $2j$, are consecutive. This labeling is not $k$-super graceful, a contradiction. Thus $k+3n+1-2j$ must be a vertex label.

\ms\nt Suppose $f(v_2) = k + 3n + 1 - 2j$, then the edges joining $v_2$ to the first $j$ vertices in $A$ have labels $2j, \ldots, j + 1$.  Now consider the next label $k + 3n - 2j > j = f(v_1)$. Since it is the largest undetermined label, it must be the label of $f(u_{j+1})$ or $f(u_{j+1}v_2)$  if $u_{j+1}$ exists. In either case, $1$ is a label, a contradiction. Thus, $u_{j+1}$ does not exist and hence $j=n$. Now $\mb{[n, 2n+1]}\cup \mb{[k+n+2, k+3n+1]}=\mb{[k, k+3n+1]}$. Hence $n=k$.

\ms\nt Suppose $f(u_{j+1})=k+3n+1-2j$ if $j<n$. Now we consider the largest undetermined label $k+3n-2j$. Since $f(u_{j+1}v_2)\ge k>2$, $f(v_2)\le k+3n-2j-2$. If $k+3n-2j=f(u_iv_2)$ for some $i$, then $i=1$ and hence $f(v_2)=2j+1$. But $f(u_jv_2)=k+3n+1-3j=f(u_{j+1}v_1)$. This is a contradiction. So $k+3n-2j$ must be a vertex label and hence $f(u_{j+2})=k+3n-2j$. By a similar argument we can show that $k + 3n+1 - 2j, \ldots, k + 3n + 2 - 3j$ are labels of $u_{j+1}, \dots u_{2j}$ and  $k + 3n +1- 3j, \ldots, k + 3n + 2 - 4j$ are labels of $u_{j+1}v_1, \dots u_{2j}v_1$, respectively. This argument can be repeated until all vertices in $A$ are labeled. Hence $n=rj$ for some $r\ge 2$. In this case, integers in $\mb{[k+rj+2, k+3rj+1]}$ are assigned to $A\cup E_1$. Now integers in $\mb{[k, k+rj+1]}\setminus \{j, f(v_2)\}$ must be assigned to $E_2$. But labels assigned to $E_2$ form $r$ subintervals with $r-1$ gaps of length $j\ge 3$. So it is impossible.
\end{proof}

\begin{lemma} Let $m, n, k \ge 2$.  If $K(m, n)$ is $k$-super graceful, and the $n$ greatest integers in $\mb{[k, k + m + n + mn - 1]}$ are labels of the partite set with $n$ vertices, then $k = n$. \end{lemma}

\begin{proof}  Let the two partite sets of $K(m, n)$ be $A = \{u_1, \ldots, u_n\}$ and $B = \{v_1, \ldots, v_m\}$.  Without loss of generality, assume that $f(v_1) < f(v_2) < \cdots < f(v_m)$, where $f$ denotes the labeling.  Furthermore, let $E_i$ be the set of edges joining $v_i$ to the vertices in $A$, for $i = 1, \ldots, m$.

\ms\nt By Theorem~\ref{thm-NAV}, $n \ge k \ge 2$.  We may assume that $f(u_j) = k + m + n + mn - j$, for $j = 1, \ldots, n$. Consider the greatest undetermined label, $k + m + mn - 1$.  If it labels a vertex in $B$, then the edge joining it to $u_n$ has label $1 < k$, which is impossible. Thus it must be the greatest edge label.  This gives $f(u_1v_1) = k + m + mn - 1$, $f(v_1) = n$, and $f(u_jv_1) = k + m + mn - j$, for $j = 1, \ldots, n$.

\ms\nt The labels of the vertices in $A$ form a block of $n$ consecutive integers consisting of the greatest $n$ labels, and the labels of the edges in $E_1$ form another block of $n$ consecutive integers consisting of the second greatest $n$ labels.  Similarly, for each $i = 2, \ldots, m$, the labels of the edges in $E_i$ form a block of $n$ consecutive integers.  These blocks of consecutive integers cannot overlap, because they consist of labels of distinct edges.  By the assumption on the labels of the vertices in $B$, the edge labels in each succeeding block must be less than the edge labels in each preceding block.  They give at most $m$ gaps, between the labels of the edges in $E_i$ and the labels of the edges in $E_{i+1}$, for $i = 1, \ldots, m - 1$, and the integers in $\mb{[k, k + m + n + mn - 1]}$ that are less than the smallest label of the edges in $E_m$.

\ms\nt Each pair of $f(v_i)$ and $f(v_{i+1})$ cannot represent consecutive integers, because otherwise $f(u_1v_{i+1}) = k + m + n + mn - 1 - f(v_{i+1}) = k + m + n + mn - 2 - f(v_i) = f(u_2v_i)$, resulting in two edges having the same label.

\ms\nt Since there are $m$ vertices in $B$, no two of which with consecutive labels, and there are at most $m$ gaps, we must have exactly $m$ gaps of one integer to label the vertices in $B$.  This forces the smallest of these vertex labels, namely $f(v_1)$, to be $k$.  As we have established that $f(v_1) = n$, this gives the desired result of $k = n$. \end{proof}

\begin{conjecture} For $m=1$, and a prime $n$ or for $n\ge m\ge 2$, $K(m,n)$ is $k$-super graceful if and only if $k=1,m,n$. \end{conjecture}

\begin{theorem} The bi-regular graph $G=K(1,\ldots,1,2)$ of order $r+2$ $(\ge 4)$ is $k$-super graceful if and only if $r=2,k=1$. \end{theorem}

\begin{proof} Let $G=K(1,\ldots,1,2)$. By~\cite[Theorem 4.8]{Lau+Shiu+Ng-kSG}, the sufficiency holds. Since the independence number of $G$ is $2$, by Theorem~\ref{thm-NAV}, $1\le k \le 2$.

\ms\nt
Suppose $G$ has a 2-super graceful labeling $f$. Let $u$ and $v$ be the vertices with degree $r$ and $w_1, \ldots, w_r$ be the vertices with degree $r+1$. For convenience we may assume $f(u) > f(v)$ and $f(w_1)>f(w_2)>\cdots >f(w_r)$. By Theorem~\ref{thm-NAV}, we have $f(u)=m$, $f(v)=m-1$, where $m$ is the largest possible label.
It follows that $m-2$ must be a label of the edge $uw_1$. Hence, $f(w_1)=2$ and $f(w_1v) = m-3$. Consequently, $m-4$ must be a label of the edge $uw_2$. Hence, $f(w_2) = 4$ which implies that $f(w_1w_2) = 2 = f(w_1)$, a contradiction.

\ms\nt Suppose $G$ has a super graceful labeling $f$ and $r\ge 3$. Again, we assume that $f(u) > f(v)$ and $f(w_1)>f(w_2)>\cdots >f(w_r)$. Note that $G$ is $K_{r+2}$ with an edge deleted. The available label set is $\mb{[1,m]}$ with  $m=(r+1)(r+4)/2\ge 14$. We consider 2 cases: (1) $f(u)=m$ and (2) $f(w_1)=m$.

\ms\nt {\bf Case 1:} $f(u)=m$. There are 3 subcases to consider: (1.1) $f(v)=m-1$, (1.2) $f(uw_1)=m-1$ and (1.3) $f(w_1)=m-1$.

\ms\nt {\bf Subcase 1.1:} Suppose $f(v)=m-1$. If $f(w_1)=m-2$, then $f(uw_1)=2, f(vw_1)=1$. This forces $f(uw_2)=m-3, f(w_2)=3, f(vw_2)=m-4, f(w_1w_2)=m-5$. This in turn forces $f(w_3)=m-6, f(uw_3)=6, f(vw_3)=5, f(w_1w_3)=4, f(w_2w_3)=m-9$. Now, $m-7$ cannot be placed, a contradiction. Therefore $m-2$ must be the label of an edge. If $f(uw_1)=m-2$, then $f(w_1)=2, f(vw_1)=m-3$. This forces $f(w_2)=m-4, f(uw_2)=4, f(vw_2)=3, f(w_1w_2)=m-6$. Now, $m-5$ cannot be placed, a contradiction.

\ms\nt {\bf Subcase 1.2:} Suppose $f(uw_1)=m-1$ so that $f(w_1)=1$. Clearly, $m-2$ must be a vertex label. If $f(v)=m-2$, then $f(vw_1)=m-3$. This forces $f(uw_2)=m-4$ or $f(w_2)=m-4$. Assume $f(uw_2)=m-4$ so that $f(w_2)=4$, $f(w_1w_2)=3$, $f(vw_2)=m-6$. However, $m-5$ cannot be placed, a contradiction. Assume $f(w_2)=m-4$ so that $f(uw_2)=4, f(vw_2)=2, f(w_1w_2)=m-5$. Now, if $m-6$ is a vertex label, we must have $f(w_3)=m-6$ and $f(w_2w_3)=f(vw_2)=2$, a contradiction. So, $m-6$ must label an edge incident to $u$ or to $v$. Suppose $f(uw_3)=m-6$ so that $f(w_3)=6, f(vw_3)=m-8, f(w_1w_3)=5, f(w_2w_3)=m-10$. Now, $m-7$ cannot be placed, a contradiction. Suppose $f(vw_3)=m-6$ so that $f(w_3)=6, f(uw_3)=m-7, f(w_1w_3)=5, f(w_2w_3)=m-10$. Now, $m-7$ cannot be placed, a contradiction.

\ms\nt If $f(w_2)=m-2$, then $f(uw_2)=2, f(w_1w_2)=m-3$. This forces $f(uw_3)=m-4, f(w_3)=4, f(w_1w_4)=3, f(w_2w_3)=m-6$. However, $m-5$ cannot be placed, a contradiction.

\ms\nt {\bf Subcase 1.3:} Suppose $f(w_1)=m-1$ so that $f(uw_1)=1$. This forces $f(uw_2)=m-2, f(w_2)=2, f(w_1,w_2)=m-3$. Hence, $f(v)=m-4$ or $f(w_3)=m-4$. Assume $f(v)=m-4$ so that $f(vw_1)=3, f(vw_2)=m-6$. However, $m-5$ cannot be placed, a contradiction. Assume $f(w_3)=m-4$ so that $f(uw_3)=4, f(w_1w_3)=3, f(w_2w_3)=m-6$. However, $m-5$ cannot be placed, a contradiction.

\ms\nt {\bf Case 2:} $f(w_1)=m$. There are 4 subcases to consider: (2.1) $f(u)=m-1$, (2.2) $f(vw_1)=m-1$, (2.3) $f(w_2)=m-1$ and (2.4) $f(w_1w_2)=m-1$.

\ms\nt {\bf Subcase 2.1:} Suppose $f(u)=m-1$. We have $f(uw_1)=1$. Clearly, we must have $m-2$ labeling (i) edge $w_1w_2$ or (ii) $vw_1$ or (iii) vertex $v$.
\begin{enumerate}[(i)]
\item If $f(w_1w_2)=m-2$, then $f(w_2)=2$ and $f(uw_2)=m-3$. This forces $m-4$ to be a vertex label. If $f(w_3)=m-4$, then there is no place to label $m-5$. If $f(v)=m-4$, then $f(vw_1)=4, f(vw_2)=m-6$. However, $m-5$ cannot be placed, a contradiction.

\item If $f(vw_1)=m-2$, then $f(v)=2$. However, $m-3$ cannot be placed, a contradiction.

\item If $f(v)=m-2$, then $f(vw_1)=2$. We must now have $m-3$ label an edge incident to $w_1$. Assume $f(w_1w_2)=m-3$, then $f(w_2)=3$, $f(uw_2)=m-4$, $f(vw_2)=m-5$. Now, $m-6$ must be the label of a vertex, say $w_3$. However, $m-7$ cannot be placed, a contradiction.
\end{enumerate}
\nt {\bf Subcase 2.2:} Suppose $f(vw_1)=m-1$ so that $f(v)=1$. Now, $m-2$ must label (i) edge $uw_1$, (ii) vertex $u$ or (iii) $w_2$.
\begin{enumerate}[(i)]
\item If $f(uw_1)=m-2$, then $f(u)=2$. We must have $f(w_2)=m-3$ so that $f(w_1w_2)=3$, $f(vw_2)=m-4$, $f(uw_2)=m-5$. This forces $f(w_1w_3)=m-6$. Hence $f(w_3)=6$, $f(vw_3)=5$, $f(uw_3)=4$. However, $m-7$ cannot be placed, a contradiction.

\item If $f(u)=m-2$, then $m-3$ cannot be placed, a contradiction.

\item If $f(w_2)=m-2$, then $f(w_1w_2)=2, f(vw_2)=m-3$. This forces $f(w_1w_3)=m-4$, $f(w_3)=4$, $f(vw_3)=3$, $f(w_2w_3)=m-6$. However, $m-5$ cannot be placed, a contradiction.
\end{enumerate}
\nt {\bf Subcase 2.3:} Suppose $f(w_2)=m-1$. So, $f(w_1w_2))=1$. This forces $f(vw_1)=m-2$ or $f(w_1w_3)=m-2$. If $f(vw_1)=m-2$, then $f(v)=2, f(vw_2)=m-3$. We consider the following 6 subcases.
\begin{enumerate}[(a)]
\item $f(uw_1)=m-4$. So, $f(u)=4, f(uw_2)=m-5$ and $m-6$ must be a vertex label, say $f(w_3)=m-6$. Now, $m-7$ cannot be placed, a contradiction.

\item $f(u)=m-4$. So, $f(uw_1)=4$, $f(uw_2)=3$. Now, $m-5$ cannot be placed, a contradiction.

\item $f(w_1w_3)=m-4$. So, $f(w_3)=4$, $f(vw_3)=f(v)=2$, a contradiction.

\item $f(w_3)=m-4$. So, $f(w_2w_3)=3$, $f(w_1w_3)=4$. Now, $m-5$ cannot be placed, a contradiction.

\item  $f(uw_2)=m-4$. So, $f(u)=3$, $f(uw_1)=f(vw_2)=m-3$, a contradiction.

\item  $f(w_2w_3)=m-4$. So, $f(w_3)=3$, $f(vw_3)=f(w_1w_2)=1$, a contradiction.
\end{enumerate}

\ms\nt If $f(w_1w_3)=m-2$, then $f(w_3)=2, f(w_2w_3)=m-3$. Now, $m-4$ can can be the label of a vertex, or of an edge incident to $w_1$ or $w_2$. We consider the following 6 subcases.

\begin{enumerate}[(a)]
\item $f(uw_1)=m-4$. So, $f(u)=4, f(uw_3)=f(w_3)=2$, a contradiction.

\item  $f(u)=m-4$. So, $m-5$ cannot be placed, a contradiction.

\item  $f(w_1w_4)=m-4$. So, $f(w_4)=4, f(w_3w_4)=f(w_3)=2$, a contradiction.

\item  $f(w_4)=m-4$. So, $f(w_1w_4)=4, f(w_2w_4)=3, f(w_3w_4)=m-6$. Now, $m-5$ cannot be placed, a contradiction.

\item $f(uw_2)=m-4$. So, $f(u)=3, f(uw_3)=f(w_1w_2)=1$, a contradiction.

\item $f(w_2w_4)=m-4$. So, $f(w_4)=3, f(w_3w_4)=f(w_1w_2)=1$, a contradiction.
\end{enumerate}

\nt {\bf Subcase 2.4:} Suppose $f(w_1w_2)=m-1$. We have $f(w_2)=1$. This forces (i) $f(u)=m-2$ or (ii) $f(w_3)=m-2$.
\begin{enumerate}[(i)]
\item If $f(u)=m-2$, then $f(uw_1)=2, f(uw_2)=m-3$. We consider the following 4 subcases.
\begin{enumerate}[(a)]
\item $f(vw_1)=m-4$. So, $m-5$ cannot be placed, a contradiction.

\item $f(v)=m-4$. So, $f(vw_1)=4$, $f(vw_2)=m-5$. This forces $f(w_1w_3)=m-6$, $f(w_3)=6$, $f(uw_3)=m-8$, $f(vw_3)=m-10$. Now, $m-7$ cannot be placed, a contradiction.
\item $f(w_1w_3)=m-4$. So, $f(uw_3)=m-6$, $f(w_2w_3)=3$. Now, $m-5$ cannot be placed, a contradiction.

\item $f(w_3)=m-4$. So, $f(uw_3)=f(uw_1)=2$, a contradiction.
\end{enumerate}

\item If $f(w_3)=m-2$, then $f(w_1w_3)=2, f(w_2w_3)=m-3$. Hence, $m-4$ must be the label of an edge incident to $w_1$. If $f(w_1w_4)=m-4$ and so $f(w_4)=4, f(w_2w_4)=3, f(w_3w_4)=m-6$. Now, $m-5$ cannot be placed, a contradiction. If $f(vw_4)=m-4$, we also get a similar contradiction.
\end{enumerate}
\end{proof}

\nt In~\cite[Theorem 4.6]{Lau+Shiu+Ng-kSG}, the authors showed that $K(1,m,r)$ is super graceful for  $m,r\ge 1$.

\begin{conjecture} For $t\ge 3$, a complete $t$-partite graph $G$ is $k$-super graceful if and only if $G=K(1,m,r)$ and $k=1$.  \end{conjecture}

\ms\nt Note that for $n\ge 3, m\ge 1$, $C_n + mK_1$ is super graceful. Particularly, $C_n + K_1$ is a super graceful $(n+1,n)$-graph. However, one can show that $C_3+P_2$ and $C_3+P_3$ are not super graceful. We now give two families of super graceful $(n+1,n)$-graphs with all odd vertex labels.

\begin{theorem} The graph $C_n + P_m$ is super graceful with all odd vertex labels if
\begin{enumerate}[(i)]
  \item $n=3$, $m\ge 4$,
  \item $n=4$, $m\ge 3$.
\end{enumerate}
\end{theorem}

\begin{proof} Consider $C_4+P_3$. Label the vertices of $C_4$ by $1,11,5,13$ consecutively and the corresponding edges by $10,6,8,12$. Label the vertices of $P_3$ by $3,7,9$ consecutively and the corresponding edges by $4,2$. Now assume $m=4r+s$, where $r\ge 1$ and $0\le s\le 3$. Let $P_m=u_1v_1u_2v_2\ldots u_{m/2}v_{m/2}$ if $m$ is even and $P_m=u_1v_1u_2v_2\ldots u_{(m-1)/2}v_{(m-1)/2}u_{(m+1)/2}$ if $m$ is odd. 

\ms\nt Suppose $n=3$.
We shall define a super-graceful labeling $f:V(C_3+P_{4r+s})\cup E(C_3+P_{4r+s})\to \mb{[1,8r+2s+5]}$ for each case.
\begin{enumerate}[(A)]
\item $s=0$.
\begin{enumerate}[(1)]
  \item Label the vertices of $C_3$ by $4r+1, 4r+3, 4r+7$ and the corresponding edges by $2,4,6$.
  \item $f(u_i)=2i+1$ and $f(u_{i+1})=f(u_i)-2$ for odd $i$, $1\le i\le 2r-1$.
  \item $f(v_{2r})=4r+5$; $f(v_1)=8r+5$; $f(v_i)=8r+9-2i$ and $f(v_{i-1}) = f(v_i)-2$ for odd $i$, $3\le i\le 2r-1$.
  \item $f(u_1v_1)=8r+2$; $f(u_iv_i)=8r+8-4i$ for $2\le i\le 2r$.
  \item $f(v_iu_{i+1})=f(u_iv_i)+2$ for odd $i$, $1\le i\le 2r-1$; $f(v_{i-1}u_i) = f(u_iv_i)-2$ for odd $i$, $3\le i\le 2r-1$.
\end{enumerate}

\item  $s=1$.
\begin{enumerate}[(1)]
  \item Label the vertices of $C_3$ by $4r+1, 4r+5, 4r+7$ and the corresponding edges by $4,2,6$.
  \item $f(u_{2r+1})=4r+3$, $f(u_i)=2i+1$ and $f(u_{i+1})=f(u_i)-2$ for odd $i$, $1\le i\le 2r-1$.
  \item $f(v_i)=8r+7-2i$ and $f(v_{i+1})=f(v_i)+2$ for odd $i$, $1\le i\le 2r-1$.
  \item $f(v_iu_{i+1})=8r+8-4i$ for $1\le i\le 2r$.
  \item $f(u_iv_i)=8r+6-4i$ and $f(u_{i+1}v_{i+1})=f(u_iv_i)+4$ for odd $i$, $1\le i\le 2r-1$.
\end{enumerate}

\item $s=2$.
\ms\nt For $m=6$, label the vertices of $C_3$ by $5,7,11$ and the corresponding edges by $2,4,6$. Label the edges of $P_6$ by $9,17,1,15,3,13$ consecutively and the corresponding edges by $8,16,14,12,10$.  Now, we consider $r\ge 2$.
\begin{enumerate}[(1)]
  \item Label the vertices of $C_3$ by $4r+3, 4r+5, 4r+9$ and the corresponding edges by $2,4,6$.
  \item $f(u_{2r-1})=4r-3$, $f(u_{2r})=4r+1$, $f(u_{2r+1})=4r-1$.
  \item $f(u_i)=2i+1$ and $f(u_{i+1})=f(u_i)-2$ for odd $i$, $1\le i\le 2r-3$.
  \item $f(v_1)=8r+9$, $f(v_{2r+1})=4r+7$, $f(v_{2r})=4r+11$.
  \item $f(v_i)=8r+9-2i$ and $f(v_{i+1})=f(v_i)+2$ for even $i$, $2\le i\le 2r-2$.
  \item $f(u_1v_1)=8r+4$, $f(v_1u_2)=8r+6$, $f(u_2v_2)=8r+2$, $f(v_{2r-2}u_{2r-1})=16$,  $f(u_{2r-1}v_{2r-1})=18$, $f(v_{2r-1}u_{2r})=14$, $f(u_{2r}v_{2r})=10$, $f(v_{2r}u_{2r+1})=12$, $f(u_{2r+1}v_{2r+1})=8$.
  \item $f(u_iv_i)=8r+10-4i$ for $3\le i\le 2r-2$.
  \item $f(v_{i-1}u_i) = f(u_iv_i) - 2$ for odd $i$, $3\le i\le 2r-3$ and $f(v_{i-1}u_i) = f(u_iv_i)+6$ for even $i$, $4\le i\le 2r-2$.
\end{enumerate}

\item $s=3$.
\begin{enumerate}[(1)]
  \item Label the vertices of $C_3$ by $4r+3, 4r+7, 4r+9$ and the corresponding edges by $4,2,6$.
  \item $f(u_1)=3, f(u_2)=5, f(u_3)=1, f(u_4)=9, f(v_1)=8r+7, f(v_2)=8r+11, f(v_3)=8r+9$.
  \item $f(u_i)=f(u_{i-2})+4$ and $f(u_{i-1})=f(u_{i-2})-2$ for even $i$, $6\le i\le 2r+2$.
  \item $f(v_i)=f(v_{i-2})-4$ and $f(v_{i-1})=f(v_i)-2$ for odd $i$, $5\le i\le 2r+1$.
  \item $f(u_1v_1)=8r+4$, $f(v_1u_2)=8r+2$, $f(u_2v_2)=8r+6$, $f(v_2u_3)=8r+10$, $f(u_3v_3)=8r+8$, $f(v_3u_4)=8r$.
  \item $f(v_{i-1}u_{i})=8r+16-4i$ for $5\le i\le 2r+2$.
  \item $f(u_iv_i)=8r+10-4i$ and $f(u_{i+1}v_{i+1})=f(u_iv_i)+4$ for even $i$, $4\le i\le 2r$.
\end{enumerate}
\end{enumerate}
\nt Suppose $n=4$. We shall define a super-graceful labeling $f:V(C_4+P_{4r+s})\cup E(C_4+P_{4r+s})\to \mb{[1,8r+2s+7]}$ for each case.

 \begin{enumerate}[(A)]
\item $s=0$.
For $m=4$, label the vertices of $C_4$ by 1, 13, 5, 15 consecutively and the corresponding edges by 12, 8, 10, 14; and label the vertices of $P_4$ by 5, 7, 3, 9 consecutively and the corresponding edges by 2, 4, 6.  Now suppose $r\ge 2$.
\begin{enumerate}[(1)]
  \item Label the vertices of $C_4$ by $1$, $8r+5$, $5$, $8r+7$ consecutively and the corresponding edges by $8r+4$, $8r$, $8r+2$, $8r+6$. 
  \item $f(u_1)=2, f(u_2)=7$, $f(u_i)=5+2i$ and $f(u_{i+1})=f(u_i)-2$ for odd $i$, $3\le i\le 2r-1$.
  \item $f(v_i)=8r+3-2i$ and $f(u_{i+1})=f(u_i)+2$ for odd $i$, $1\le i\le 2r-1$.
  \item $f(u_1v_1)=8r-2$, $f(v_1u_2)=8r-6$, $f(u_2v_2)=8r-4$.
  \item $f(v_iu_{i+1})=8r+4-4i$ for $4\le i\le 2r$.
  \item $f(u_{i}v_{i})=f(v_iu_{i+1})-2$ for odd $i$, $3\le i\le 2r-1$.
  \item $f(u_iv_i)=f(v_{i-1}u_{i})+2$ even $i$, $4\le i\le 2r$.
\end{enumerate}

\item $s=1$.
\begin{enumerate}[(1)]
  \item Label the vertices of $C_4$ by $1$, $8r+5$, $3$, $8r+9$ consecutively and the corresponding edges by $8r+4$, $8r+2$, $8r+6$, $8r+8$.
  \item $f(u_i)=8r+9-2i$ and $f(u_{i+1})=f(u_i)-6$ for odd $i$, $1\le i\le 2r+1$.
  \item $f(v_i)=5+2i$ and $f(v_{i+1})=f(v_i)-2$ for odd $i$, $1\le i\le 2r-1$.
  \item $f(u_iv_i)=8r+4-4i$ for $1\le i\le 2r$.
  \item $f(v_iu_{i+1})=f(u_iv_i)-6$ for odd $i$, $1\le i\le 2r-1$.
  \item $f(v_{i-1}u_i)=f(u_{i-1}v_{i-1})+2$ for odd $i$, $3\le i\le 2r+1$.
\end{enumerate}

\item $s=2$.
\begin{enumerate}[(1)]
  \item Label the vertices of $C_4$ by $1$, $8r+7$, $3$, $8r+11$ consecutively and the corresponding edges by $8r+6$, $8r+4$, $8r+8$, $8r+10$.
  \item $f(u_1)=8r+9$, $f(u_2)=8r+5$, $f(u_3)=8r+1$, $f(u_i)=f(u_{i-2})-4$ and $f(u_{i-1})=f(u_i)+6$ for odd $i$, $5\le i\le 2r+1$.
  \item $f(v_{m/2})=4r+7$, $f(v_i)=5+2i$ and $f(v_{i+1})=f(v_i)-2$ for odd $i$, $1\le i\le 2r-1$.
  \item $f(u_1v_1)=8r+2$, $f(v_1u_2)=8r-2$, $f(u_2v_2)=8r$ and $f(u_{i}v_{i-1})=8r+8-4i$ for $3\le i\le 2r+1$.
  \item $f(u_{2r+1}v_{2r+1})=2$, $f(u_iv_i)=f(v_{i-1}u_i)-6$ for odd $i$, $3\le i\le 2r-1$.
  \item $f(u_iv_i)=f(v_{i-1}u_i)+2$ for even $i$, $4\le i\le 2r$.
\end{enumerate}

\item $s=3$. 
\begin{enumerate}[(1)]
  \item Label the vertices of $C_4$ by $1,8r+9,3,8r+13$ consecutively and the corresponding edges by $8r+8, 8r+6, 8r+10, 8r+12$.
  \item $f(u_{2r+2})=4r+3$, $f(v_{2r+1})=4r+7$, $f(v_{2r+1}u_{2r+2})=2$.
  \item $f(u_i)=8r+13-2i$ and $f(u_{i+1})=f(u_i)-6$ for odd $i$, $1\le i\le 2r+1$.
  \item $f(v_{2r+1})=4r+7$, $f(v_i)=5+2i$ and $f(u_{i+1})=f(u_i)-2$ for odd $i$, $1\le i\le 2r-1$.
  \item $f(u_iv_i)=8r+8-4i$ for $1\le i\le 2r+1$.
  \item $f(v_{i-1}u_i)=f(u_iv_i)-2$ and $f(v_iu_{i+1})=f(u_iv_i)+2$ for even $i$, $2\le i\le 2r$.
\end{enumerate}
\end{enumerate}
\nt Clearly, each function defined above is a super graceful labeling of the respective graph having all odd vertex labels. \end{proof}

\nt By Theorem~\ref{thm-SGnG}, we have

\begin{corollary} The graphs $C_3 + P_m, m\ge 4$ and $C_4+P_m, m\ge 3$ are graceful.   \end{corollary}


%
%
%
\begin{theorem}\label{thm-2paths-1} For $n\ge 2$, the graph $P_{n+r}+P_n$ is 2-super graceful with edge-label set $\mb{[2,2n+r-1]}$ for $r=2,3$.   \end{theorem}

\begin{proof} Suppose $r=2$. Begin with $P_{2n+1} = u_1u_2\ldots u_{2n+1}$. First label the edges $u_1u_2, u_2u_3, \ldots, u_{2n}u_{2n+1}$ by $2, 3, \ldots, 2n+1$. Next label the vertices $u_{2n+1}, u_{2n-1}, u_{2n-3}, \ldots, u_1$ by $2n+2, 2n+3, 2n+4, \ldots, 3n+2$. Now, label the vertices $u_2, u_4, u_6, \ldots, u_{2n}$ by $3n+4, 3n+5, 3n+6, \ldots, 4n+3$. Finally, delete the edge $u_nu_{n+1}$ that has label $n+1$ and join a new vertex $u$ to $u_{2n+1}$ by labeling $u$ and $uu_{2n+1}$ by $3n+3$ and $n+1$ respectively.

\ms\nt For $r=3$, begin with $P_{2n+2}$, we can also get a 2-super graceful labeling for $P_{n+3} + P_n$ in a similar way.
\end{proof}

\begin{theorem}\label{thm-2paths-2} For $n\ge 2$, (i) $P_{n+2} + P_n$ is $(n+1)$-super graceful with edge-label set $\mb{[2n+2,4n+1]}$ and (ii) $P_{n+3}+P_n$ is $(n+1)$-super graceful with edge-label set $\mb{[2n+3,4n+3]}$.   \end{theorem}

\begin{proof} (i) Begin with $P_{2n+1} = u_1u_2\ldots u_{2n+1}$. First label the vertices $u_2, u_4, u_6, \ldots, u_{2n}$ by $n+1, n+2, n+3, \ldots, 2n$. Next label the edges $u_{2n+1}u_{2n}, u_{2n}u_{2n-1}, \ldots, u_2u_1$ by $2n+2, 2n+3, \ldots, 4n+1$. Now, label the vertices $u_{2n+1}, u_{2n-1}, u_{2n-3}, \ldots, u_1$  by $4n+2, 4n+3, 4n+4, \ldots, 5n+2$. Finally, delete the edge $u_{n+1}u_{n+2}$ that has label $3n+1$ and join a new vertex $u$ to $u_1$ by labeling $u$ and $uu_1$ by $2n+1$ and $3n+1$ respectively. 

\ms\nt (ii) Begin with $P_{2n+2}$, we can also get an $(n+1)$-super graceful labeling for $P_{n+3} + P_n$ in a similar way.
\end{proof}

\nt In~\cite[Construction C4]{Lau+Shiu+Ng-kSG}, the authors gave a way to construct infinitely many super graceful bipartite graphs. The approach can be extended to obtained infinitely many $k$-super graceful bipartite graphs.

\ms\nt{\it Approach A1.} Begin with vertices $u_i$ $(0\le i\le n)$. Choose an integer $d > k\ge 1$.
\begin{enumerate}[(a)]
  \item Label $u_i$ by $k+id$.
  \item For $1\le j\le d-1$, add a vertex $v_j$ and join it to each of $u_i$.
  \item Label edge $u_iv_j$ by $k+j+(n-i)d$ and vertex $v_j$ by $2k+j+nd$.
  \item Delete edge $u_iv_j$ if its label is also one of the vertex labels.
  \item For $r=1,2,\ldots$, introduce $d$ new vertices with labels $(r+2)k+(nr+n+r)d+s-1$, $1\le s\le d$. Join each of them to  each  $u_i$, $0\le i\le n$. The induced edge labels are $(r+1)k+(nr+n+r-i)d+s-1$.
  \item Delete each new edge in (e) if its label is one of the new vertex labels.
\end{enumerate}

\ms\nt{\it Approach A2.} Begin with integers $u_i$ $(0\le i\le n)$. Choose an integer $2\le d\le k$
\begin{enumerate}[(a)]
  \item Label $u_i$ by $k+id$.
  \item For $1\le j\le k$, add a vertex $v_j$ and join it to each of $u_i$.
  \item Label edge $u_iv_j$ by $k+j+(n-i)d$ and vertex $v_j$ by $2k+j+nd$.
  \item In Step (c), if an edge label is also a vertex label, delete the corresponding edge. If a label is assigned to more than 1 edge, delete all but one of the edges.
  \item For $r=1,2,\ldots$, introduce $k$ new vertices with labels $2(r+1)k+(r+1)d+s$, $1\le s\le k$. Join each of them to  each $u_i$, $0\le i\le n$. The induced edge labels are $(2r+1)k+rd+s$.
  \item In Step (e), if a label is assigned to more than 1 edge, delete all but one of the edges.
\end{enumerate}

\nt One can verify that the bipartite graphs we have thus obtained are $k$-super graceful. Moreover, adding and assigning appropriate labels to the edges $v_1v_2, v_1v_3, \ldots, v_1v_{d-1}$ in Approach A1, and to the edges $v_1v_2, v_1v_3, \ldots, v_1v_k$ in Approach A2 give us $1$-, $2$-, \ldots, $(k-1)$-super graceful tripartite graphs respectively.

\ms\nt For $k\ge 3$, we also obtained three other $k$-Skolem sequences of length $2k-1$. Hence, we end this paper with the following problem.

\begin{problem} Determine the total number of distinct $k$-Skolem sequence of length $n$ for all possible $n,k$.  
\end{problem}

\ms\nt {\bf Acknowledgement} Part of this paper was done during the first author's visit to Harbin Engineering University in 2017. He is grateful to the university for providing full financial support and to Universiti Teknologi MARA for granting the leaves.

\bibliographystyle{plainnat}

\end{document}